\documentclass[a4paper,11pt]{article}%
\usepackage{amssymb}
\usepackage{epsfig}
\usepackage{graphicx}
\usepackage{amsmath}
\usepackage{amsfonts}%
\begin{document}

\title{\textbf{Positive forward rates in the maximum smoothness framework}}
\author{Juli\'{a}n Manzano\thanks{corresponding author: Juli\'{a}n Manzano, Division
of Applied Physics. Institute of Physics and Measurement
Technology. Link\"{o}ping University, S-581 83 Link\"{o}ping,
Sweden; e-mail:
manzano@ifm.liu.se; phone: +46.737.353.373.}\\Division of Applied Physics.\\Institute of Physics and Measurement Technology.\\Link\"{o}ping University.\\and \\J\"{o}rgen Blomvall \\
Division of Optimization. Department of
Mathematics.\\Link\"{o}ping University.}
\date{}
\maketitle

\begin{abstract}
In this article we present a non-linear dynamic programming algorithm for the
computation of forward rates within the maximum smoothness framework. The
algorithm implements the forward rate positivity constraint for a
one-parametric family of smoothness measures and it handles price spreads in
the constraining dataset. We investigate the outcome of the algorithm using
the Swedish Bond market showing examples where the absence of the positive
constraint leads to negative interest rates. Furthermore we investigate the
predictive accuracy of the algorithm as we move along the family of smoothness
measures. Among other things we observe that the inclusion of spreads not only
improves the smoothness of forward curves but also significantly reduces the
predictive error.

\end{abstract}

\vfill\vbox{ May 2003\null\par }

\newpage

\section{Introduction}

\label{intro}Within the financial industry forward rate curves play a central
role in fixed-income derivative pricing and risk management. Notwithstanding,
such curves are not empirical or directly measurable objects but rather useful
abstract concepts from where observed prices can be derived. Moreover, given a
finite set of market prices, we can construct in general an infinite number of
compatible forward rate curves\footnote{We assume here that prices do not
allow for arbitrage. If this is not the case the compatibility condition has
to be relaxed.}. To avoid such ambiguity several approaches have been proposed
in the literature trying to capture a \textquotedblleft
reasonable\textquotedblright\ or \textquotedblleft natural\textquotedblright%
\ functional form within the set of compatible possibilities.

In historical order the first kind of methods proposed to solve this problem
make use of the so-called parametric approach. In this approach a particular
functional form for the forward rate curve is assumed leaving a certain number
of free parameters to be fixed from the calculation of a given set of quoted
prices. An extensive literature exists advocating for this approach. We can
cite as examples the works of McCulloc \cite{McCulloch}, Vasicek and Fong
\cite{Vasicek}, Chambers, Carleton and Waldman \cite{Chambers}, Shea
\cite{Shea}, Nelson and Siegel \cite{Nelson} and more recently the works of
Svensson \cite{Svensson}, Fisher, Nychka and Zervos \cite{Fisher} and Waggoner
\cite{Waggoner}. In most of these works we notice the privileged role played
by polynomial and exponential splines as the preferred functional forms for
the forward rate curves.

The second kind of methods has been termed in the literature as non-parametric
or maximum-smoothness approach. Here instead of advocating for an a priori
functional form for the forward rate a given measure of smoothness is chosen
and then the forward rate curve is obtained as the one maximizing this measure
subject to the constraints imposed by market prices. Examples where these
methods have been investigated include the works of Adams and Van Deventer
\cite{Adams}, Delbaen and Lorimier \cite{Delbaen}, Kian Guan Lim, Qin Xiao et
al \cite{Lim Xiao, Lim Xiao Ang}, Frishling and Yamamura \cite{Frishling} and
Yekutieli \cite{Yekutieli}. In these works three different smoothing measures
have been proposed. We also have the works of Forsgren \cite{Forsgren} and
Kwon \cite{Kwon} that generalize these methodologies and clarify the
connection between splines and certain smoothness measures. Finally we point
out the work of Wets, Bianchi and Yang \cite{Wets} that can be located
somewhere in-between both approaches since here the number of functional
parameters is finite (albeit arbitrarily large) and the functional behavior is
restricted to a subfamily of $C^{2}$ curves.

The purpose of this article is two-fold: Firstly we want to present an
efficient maximum-smoothing algorithm that handles the presence of spreads and
implements the positivity constraint. Secondly we want to investigate the
predictive power of a linear combination of two quadratic measures, namely the
one proposed by Delbaen et al. and Frishling et al. \cite{Delbaen, Frishling}
and the one by Adams and van Deventer \cite{Adams}. Here it is worth remarking
that once the compatibility with market prices is fulfilled the only guiding
principle that should be taken as definition of \textquotedblleft
reasonable\textquotedblright\ or \textquotedblleft natural\textquotedblright%
\ is the predictive power and not other ad-hoc criteria.

In this article we will only use as constraining data coupon bearing bonds.
The inclusion of treasury bills, zero coupon bonds or bill futures is
straightforward and amounts to adding the corresponding linear constraints.
Since our objective in this article is focusing on an algorithm dealing with
non-linear constraints and inequality constraints (spreads and positivity
constraint) we have not included such data.

With these objectives in mind we organize the article as follows:

In section \ref{sec algorithm} we present the objective function that we will
use throughout the article and we establish the basic notation. In this
section we present a sketch of the complete algorithm leaving the details for
the appendices. In section \ref{sec results} we present the results of the
article including examples where the absence of the positivity constraint or
the adequate spreads leads to negative rates. Here we present also a study of
the predictive power of the one-parametric family of smoothness measures that
include as extreme cases the measures used by Delbaen et al., Frishling et al.
and Adams and Van Deventer. Finally in section \ref{sec conclusions} we
present the conclusions.

\section{The algorithm}

\label{sec algorithm}A bond, $j$, is an instrument that gives future coupons,
$c_{ij}$, at time stages $R_{i}^{(j)},i=1,\ldots,n_{j}-1$ and a final payment,
$c_{n_{j}j}=N_{j}$\footnote{By final payment $N_{j}$ we mean the complete last
cash flow of bond $j$, typically that includes a principal plus a last
coupon.}. The bond price, $P_{j}$, can be determined from the discrete forward
rate curve, $f_{r},r=1,\ldots,R_{n_{j}}^{(j)}$, as follows
\begin{equation}
P_{j}=\sum\limits_{i=1}^{n_{j}}c_{ij}\exp\left(  -\sum_{r=1}^{R_{i}^{\left(
j\right)  }-1}f_{r}\xi_{r}\right)  , \label{bondPrice}%
\end{equation}
where $\xi_{r}$ is the length of the time period between time stage $r$ and
$r+1$ (in our implementation we have used $\xi_{r}=1$ day).

The objective function, or smoothing measure\footnote{Note that since we do
not include a global minus sign we have to perform a minimization and not a
maximization. With this sign, that is the one used in the literature,
\textquotedblleft rugosity\textquotedblright\ measure would be a more
appropriate name.}, is defined as a linear combination of the one used by
Delbaen et al. and Frishling et al. (DF) \cite{Delbaen, Frishling} and the one
used by Adams and Van Deventer (AD) \cite{Adams}
\begin{equation}
W:=\frac{\gamma}{2}\sum_{r=1}^{n-1}\left(  \frac{f_{r+1}-f_{r}}{\xi_{r}%
}\right)  ^{2}\xi_{r}+\frac{\varphi}{2}\sum_{r=2}^{n-1}\left(  \frac
{2}{\left(  \xi_{r-1}+\xi_{r}\right)  }\left(  \frac{f_{r+1}-f_{r}}{\xi_{r}%
}-\frac{f_{r}-f_{r-1}}{\xi_{r-1}}\right)  \right)  ^{2}\xi_{r},
\label{objective}%
\end{equation}
The first term in Eq.(\ref{objective}) (DF) is a discrete approximation of the
integral of the square of the first derivative of $f$ and the second term (AD)
is a discrete approximation of the integral of the square of the second
derivative of $f$. This objective function is to be minimized subject to the
consistency constraints%
\begin{align}
0  &  =\rho_{j}+\sum_{r=1}^{R_{n_{j}}^{\left(  j\right)  }-1}f_{r}\xi_{r}%
-\ln\left(  v_{j}\right)  ,\label{eq_const}\\
f_{r}  &  \geq0,\qquad\rho_{j}^{b}\leq\rho_{j}\leq\rho_{j}^{a},
\label{ineq_constr}%
\end{align}
where $\rho_{j}=\ln\left(  P_{j}/N_{j}\right)  $ is to be determined along
with $f$ and where we have used the definitions
\begin{align}
v_{j}:=  &  \sum\limits_{i=1}^{n_{j}}\alpha_{ij}\exp\left(  \sum
_{r=R_{i}^{\left(  j\right)  }}^{R_{n_{j}}^{\left(  j\right)  }-1}f_{r}\xi
_{r}\right)  ,\qquad\rho_{j}^{b}:=\ln\left(  \alpha_{0j}^{b}\right)
,\qquad\rho_{j}^{a}:=\ln\left(  \alpha_{0j}^{a}\right)  ,\nonumber\\
\alpha_{ij}:=  &  \frac{c_{ij}}{N_{j}},\qquad\alpha_{0j}^{b}:=\frac{P_{j}^{b}%
}{N_{j}},\qquad\alpha_{0j}^{a}:=\frac{P_{j}^{a}}{N_{j}}\qquad\alpha_{n_{j}%
j}:=1, \label{definitions}%
\end{align}
with $P_{j}^{b},$ $P_{j}^{a}$ the respective bid and ask prices of bond $j$
($j=1,\cdots,m$)\footnote{Note that this optimization problem is non-convex
and therefore several local minima may exist.}. Note that the equality
constraint given by Eq.(\ref{eq_const}) is just Eq.(\ref{bondPrice}) rewritten
taking logarithms and using the definitions (\ref{definitions}).
Eq.(\ref{ineq_constr}) introduces two inequality constraints. The first one is
the positivity constraint over the forward rate curve and the second one is
the requirement that the single price given by Eq.(\ref{bondPrice}) must lie
in-between the bid and ask prices.

We define the spread of bond $j$ as the quantity $\rho_{j}^{a}-\rho_{j}^{b}$.
We take the largest time to maturity in Eq.(\ref{objective}) equal to the
largest time to maturity in the constraining dataset, namely%
\[
n:=\max_{j}\left(  R_{n_{j}}^{\left(  j\right)  }\right)  .
\]

\begin{figure}[ptbh]
\begin{center}
\includegraphics[width=15cm]{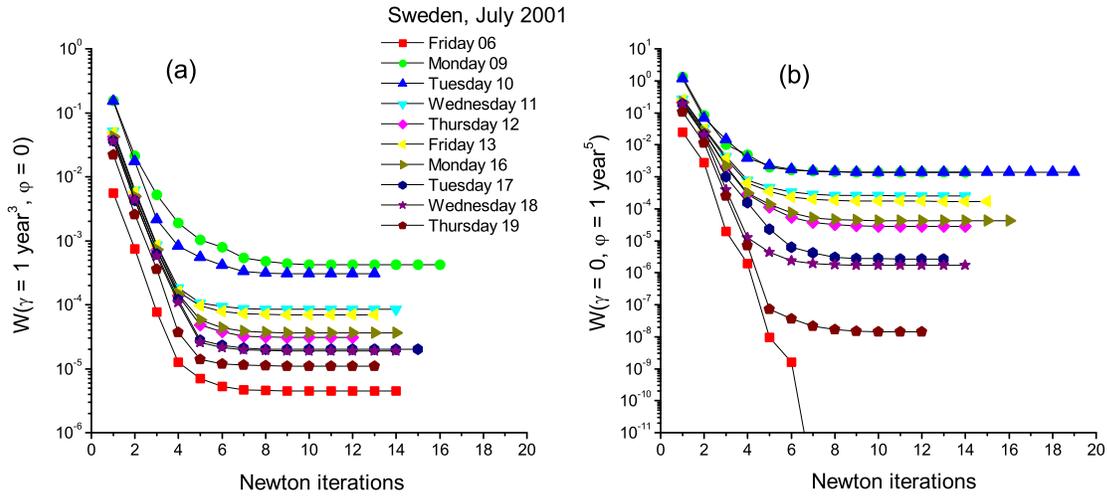}
\end{center}
\caption{Plot of the objective function $W$ (see Eq. (\ref{objective})) vs.
Newton iterations (starting from the seeds: $f^{[0]}=0.04/year$, $\rho
_{j}^{[0]}=(\rho_{j}^{a}+\rho_{j}^{b})/2$) . The plots correspond to the
forward rate curves of Fig.(\ref{Graph1}). The converge behavior is affected
by the termination criterion dictated by Eq.(\ref{criterion}) and by the
initial values of the log barriers dictated by Eq.(\ref{values}). Note the
exponential convergence of the algorithm in the first steps. Note also that on
Friday 06 the minimal solution for AD measure is, for all practical purposes,
a straight line. }%
\label{Graph2}%
\end{figure}

\begin{figure}[ptbhptbh]
\begin{center}
\includegraphics[width=15cm]{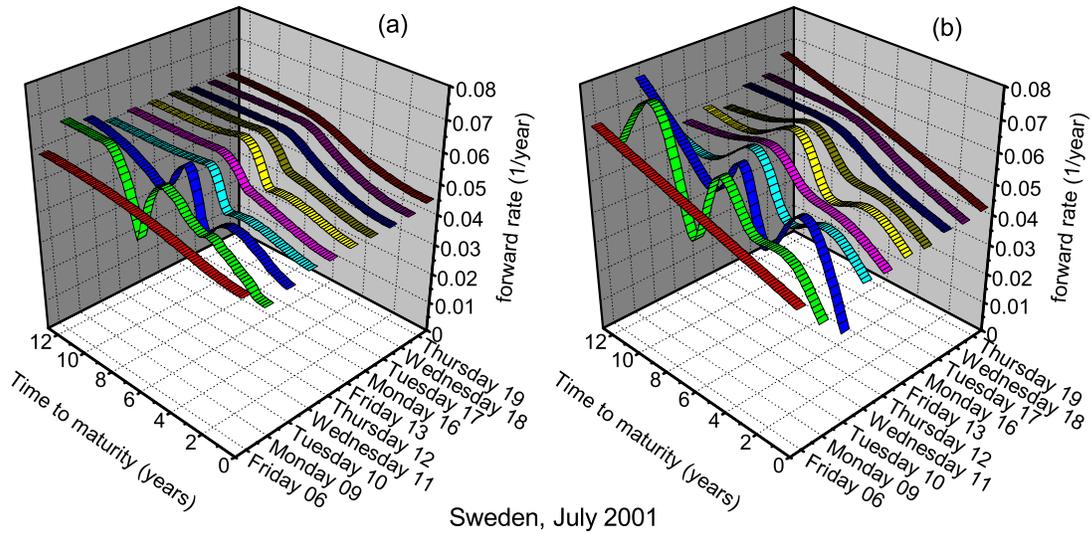}
\end{center}
\caption{Forward rate curves obtained from the prices of 11 Swedish government
bonds in a two-week period of July 2001. The curves in plot (a) are calculated
using DF measure ($\gamma=1$ year$^{3}$, $\varphi=0$) while the ones in plot
(b) correspond to the AD measure ($\gamma=0$, $\varphi=1$ year$^{5}$). In both
plots we have used a tolerance spread of 1\% ($\rho_{j}^{a}-\rho_{j}^{b}%
=0.01$). }%
\label{Graph1}%
\end{figure}

The constraints reflecting bond prices (\ref{bondPrice}) have been rewritten
in a way such that they become linear when no coupons are present ($v_{j}=1$
for a zero--coupon bond $j$). Constraints given by Eqs.(\ref{eq_const}) and
(\ref{ineq_constr}) are moved to the objective function defining%
\begin{equation}
Z:=W+\sum_{j=1}^{m}\lambda_{j}\left[  \rho_{j}+\sum_{r=1}^{R_{n_{j}}^{\left(
j\right)  }-1}f_{r}\xi_{r}-\ln\left(  v_{j}\right)  \right]  -\mu\sum
_{r=1}^{n}\ln\left(  f_{r}\right)  -\tilde{\mu}\sum_{j=1}^{m}\left(
\ln\left(  \rho_{j}-\rho_{j}^{b}\right)  +\ln\left(  \rho_{j}^{a}-\rho
_{j}\right)  \right)  , \label{total_obj}%
\end{equation}
with the Lagrange multipliers $\lambda_{j},$ $j=1,\ldots,m$ and the
logarithmic barriers with parameters $\mu>0$ and $\tilde{\mu}>0$ (in the
solution procedure we take $\mu\rightarrow0$, $\tilde{\mu}\rightarrow0$). The
use of log barriers to deal with inequality constraints is a standard
methodology in interior point methods for optimization problems \cite{Nocedal}%
. An explanation of this methodology adapted to our problem is given in
appendix A. Briefly the minimization algorithm is structured as follows:%
\begin{align}
\text{step 0}\text{: }  &  \text{Initialize log barriers with coefficients
}\mu^{\left[  0\right]  }\text{, }\tilde{\mu}^{\left[  0\right]  }\text{ and
set (}f,\rho\text{) }=\text{(}f^{\left[  0\right]  },\rho^{\left[  0\right]
}\text{) }\nonumber\\
&  \text{with the seed (}f^{\left[  0\right]  },\rho^{\left[  0\right]
}\text{) satisfying the inequality constraints (\ref{ineq_constr}). Let
}k=0.\nonumber\\
\text{step 1}\text{: }  &  \text{Make a second order approximation of }Z\text{
at (}f^{\left[  k\right]  },\rho^{\left[  k\right]  }\text{) (see appendix
A).}\nonumber\\[0.03in]
\text{step 2}\text{: }  &  \text{Determine the newton step (}\hat{f}^{\left[
k+1\right]  },\hat{\rho}^{\left[  k+1\right]  }\text{) using dynamic
programming (see appendix B).}\nonumber\\[0.03in]
\text{step 3}\text{: }  &  \text{Modify (}\hat{f}^{\left[  k+1\right]  }%
,\hat{\rho}^{\left[  k+1\right]  }\text{) to get a solution (}f^{\left[
k+1\right]  },\rho^{\left[  k+1\right]  }\text{) that satisfies the
inequality}\nonumber\\
&  \text{constraints (\ref{ineq_constr}). Update log barriers (}\mu^{\left[
k\right]  }\rightarrow0\text{ and }\tilde{\mu}^{\left[  k\right]  }%
\rightarrow0\text{ as }k\rightarrow\infty\text{) and}\nonumber\\
&  \text{calculate the values of }W\text{ and of constraints (\ref{eq_const}).
Check if a termination criterion}\nonumber\\
&  \text{is satisfied, otherwise let }k=k+1\text{ and go to step 1 (see
appendix C).} \label{algorithm}%
\end{align}
Computing times involved in step 2 are summarized in subsection B.1. The
solution typically stabilizes in approximately 6 iterations as can be seen in
Fig.(\ref{Graph2}). On a Pentium 4, 2.4 Ghz computer the algorithm coded in
C++ takes around 1/4 sec. to compute a forward rate curve like any of the ones
seen in Fig.(\ref{Graph1}). 

\section{Consistency and predictability}

\label{sec results}In this section we present some examples of the behavior of
algorithm (\ref{algorithm}) and we investigate the predictive power of the
resulting forward rate curves. In Fig.(\ref{Graph1}) we present a series of
forward rate curves calculated using DF and AD smoothing measures. There we
can see that for both measures the resulting curves share some similar traits
like the positions of most peaks and dips. Clear differences between both sets
of curves are found at their end-points and in their behavior in the presence
of high spreads. In the set of curves obtained from the DF measure we have
vanishing first derivative at the end-points and curves that tend to constants
for high spreads. For the AD measure we have vanishing second and third
derivatives at end-points \cite{Kwon} and curves for high enough spreads given
by straight lines. From the financial point of view these features are, in
principle, just different aesthetic possibilities. In order to choose a
particular measure the guiding principles should be, in the first place, the
fulfillment of the consistency constraints given by Eqs.(\ref{eq_const}%
-\ref{ineq_constr}) and after this is guaranteed the predictive performance.

\begin{figure}[ptbh]
\begin{center}
\includegraphics[width=16cm]{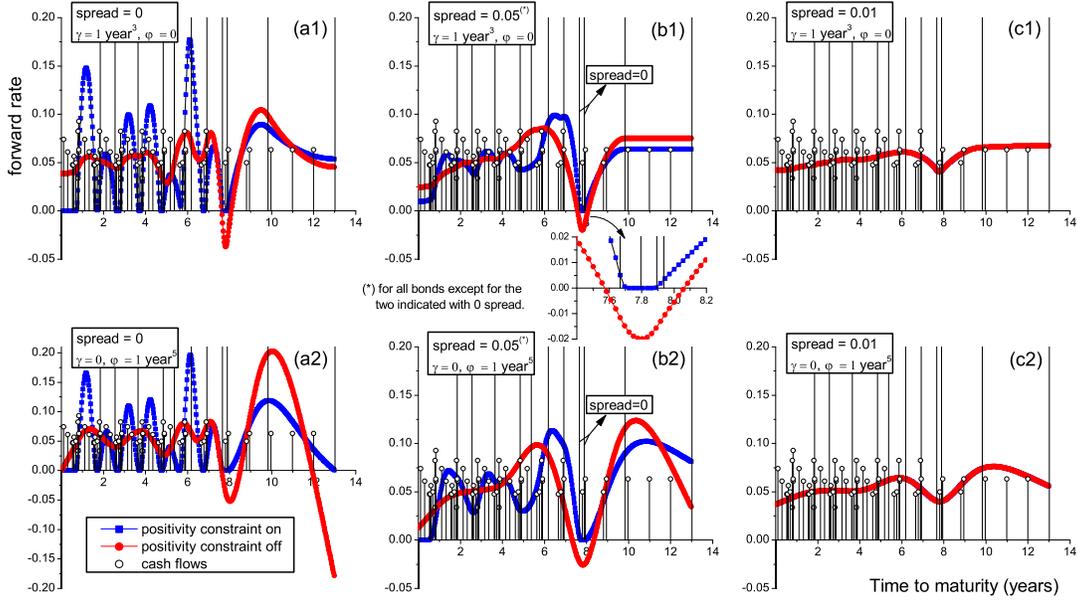}
\end{center}
\caption{Forward rate curves obtained from the prices of 11
Swedish government bonds on Monday, July 09, 2001. The vertical
lines indicate cash flows. Highest lines indicate the 11 cash
flows at maturities (normalized to 1) and the smaller ones coupon
cash flows. In these plots we show the effect of the positivity
constraint with three spread patterns both for the DF (in (a1),
(b1) and (c1)) and the AD (in (a2), (b2) and (c2)) measures. In
plots (c1) and (c2) we see that a spread of 1\% is enough to
generate positive and sensible curves for both measures (the
positivity constraint has no effect in this case). In plots (a1)
and (a2) we take null spreads observing that for both measures
large negative rates are obtained (conspicuously for AD measure).
Also in these plots we see that even though the positivity
constraint can be fulfilled the resulting curves have large
oscillations. Negative rates or large oscillations in the
positively constrained curves can be related to price patterns
that are close to violate absence of arbitrage (the precise nature
of this relation will be studied elsewhere). Plots (b1) and (b2)
explores this fact keeping null spreads in two particular bonds
(bonds SO 1043 and SO 1034 in Tables \ref{Table1} and
\ref{Table2}) and allowing for a large 5 \% in the rest (almost
unconstraining this subset). Doing so we observe in (b1) and (b2)
curves that are still negative in the region
where the nominal cash flows of the 0-spread bonds take place. }%
\label{Graph3}%
\end{figure}

Let us start now with the analysis of consistency. As can be found in
\cite{Kwon} if we do not consider the positivity constraint the local minima
of objective (\ref{objective}) are given by exponential splines with exponents
$\pm\sqrt{\gamma/\varphi}$ or polynomial splines of order 2 or 4 when
$\varphi$ or $\gamma$ are respectively zero. The main problem with these
exponential or polynomial splines is that there is no warranty that they
fulfill the positivity constraint. Negative rates are not admissible in the
absence of arbitrage opportunities and the risk of obtaining this unwanted
feature is illustrated in Fig.(\ref{Graph3}). In this figure we have
concentrated on the Swedish bond data on Monday, July 09, 2001. There we have
tested three spread patterns for both DF and AD measures with and without the
positivity constraint. From these plots it is evident that the inclusion of
spreads in the calculation of the forward rates is a necessary ingredient that
can have a major impact in the resulting functional behavior.

\begin{figure}[ptbhptbh]
\begin{center}
\includegraphics[width=15cm]{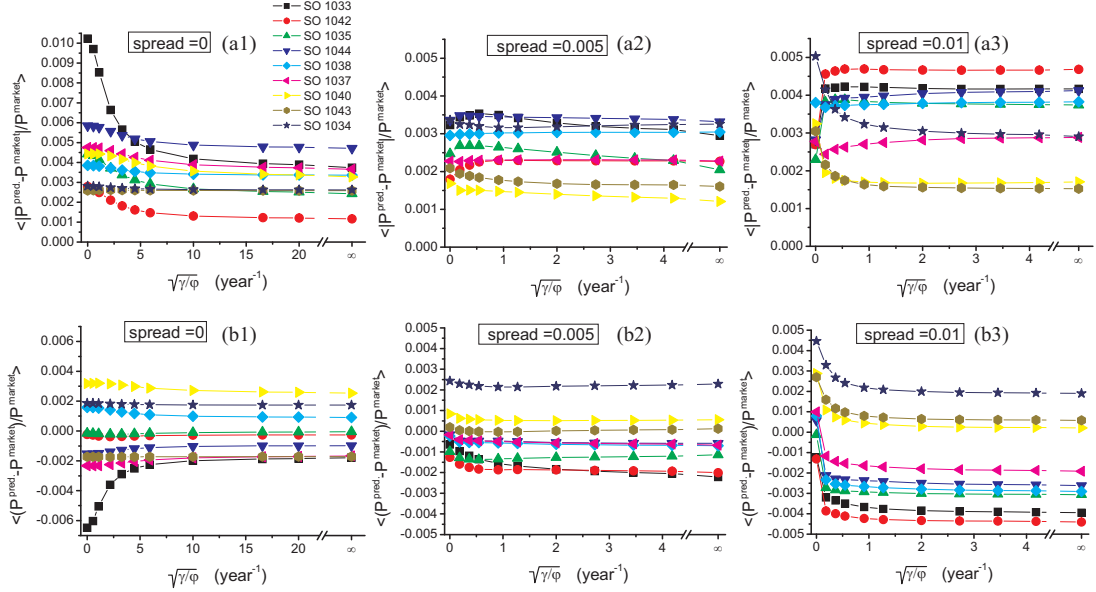}
\end{center}
\caption{Plots indicating some statistical features of price prediction errors
as a function of $\sqrt{\gamma/\varphi}$. These prediction errors are obtained
removing the indicated bonds from the constraining dataset. The statistical
sample comprise the forward curves corresponding to 335 consecutive trading
days starting on Wednesday, November 08, 2000 and ending on Thursday, March
07, 2002. Within this set of days we calculate the relative errors of the
predicted prices given by $\left(  P^{pred.}-P^{market}\right)  /P^{market}$
where $P^{market}$ is the actual price of the bond given by the market and
$P^{pred.}$ is the price predicted using the prices of all other bonds in the
dataset. Hence in plots (a1), (a2) and (a3) we show the average of the
absolute value of these errors and the plots (b1), (b2) and (b3) just their
average. Consecutive columns show averages obtained using spreads of $0\%$,
$0.5\%$ and $1\%$ in the constraining dataset. The set of predicted bonds
shown in this figure is given by the first 9 bonds of Table \ref{Table2} . The
remaining bonds of this table (the last two bonds with the larger maturities)
are analyzed in Fig.(\ref{Graph5}). Note how prediction errors tend to
decrease when spreads are considered albeit not significantly (compare with
Fig.(\ref{Graph5})). Note also that in the absence of spreads the DF measure
($\sqrt{\gamma/\varphi}=\infty$) systematically exhibits a better performance
than the AD one ($\sqrt{\gamma/\varphi}=0$)).}%
\label{Graph4}%
\end{figure}

\begin{figure}[ptbhptbhptbh]
\begin{center}
\includegraphics[width=15cm]{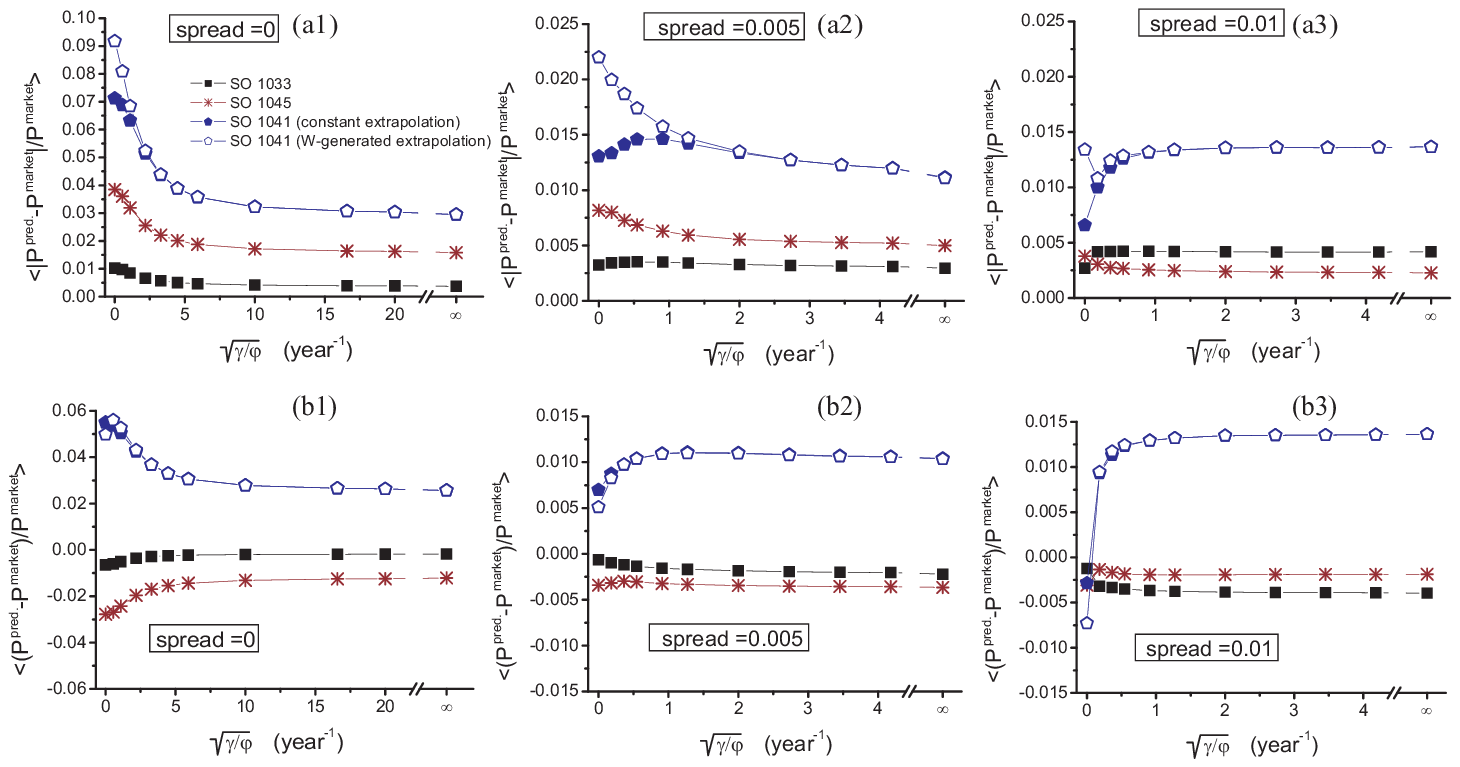}
\end{center}
\caption{In this figure we consider the same kind of statistical
features as in Fig.(\ref{Graph4}) but now the plots correspond to
the last two bonds of Table \ref{Table2} that are the ones with
larger maturities within our dataset (we have also included the
bond with the shortest maturity to facilitate the comparison with
Fig.(\ref{Graph4})). For the bond with the largest maturity the
forward curve has to be extrapolated to reach its last cash flows
(see Fig.(\ref{Graph3})). We consider two extrapolation
possibilities, one where we continue the forward rate curve from
the last cash flow as a constant and the other where the curve is
dictated by the minimization of functional $W$ even beyond the
last cash flow in the constraining dataset. Note that constant
extrapolation gives better results than W-generated extrapolation
in the region close to the AD measure ($\sqrt{\gamma/\varphi}=0$).
Again like in Fig.(\ref{Graph4}) we observe that in the absence of
spreads the DF measure ($\sqrt{\gamma/\varphi}=\infty$)
systematically exhibits a better performance than the AD one. The
most striking difference with Fig.(\ref{Graph4}) is that here we
observe that for these two long maturing bonds the introduction of
spreads drastically reduce the prediction error all along the
family of
measures.}%
\label{Graph5}%
\end{figure}

Once we have an algorithm that insures the consistency of forward rates we can
concentrate on the predictive accuracy of different measures. However, before
starting the analysis of this issue let us make a brief digression to comment
a point regarding measure (\ref{objective}). If we want to have both $\gamma$
and $\varphi$ different from zero and we want to compare the effects of each
term it is important to realize that DF and AD measures scale differently
under changes of time units. In other words, $\gamma$ and $\varphi$ have
different units. A practical way to define their units is to consider the
objective $W$ as an adimensional quantity. By doing so and remembering that
$f$ has units of inverse time, it is immediate to obtain that $\gamma$ has
units of time$^{3}$ and $\varphi$ units of time$^{5}$. The importance of
keeping this in mind becomes apparent in results like the ones presented in
Figs.(\ref{Graph4}) and (\ref{Graph5}).

Figs.(\ref{Graph4}) and (\ref{Graph5}) summarize our results regarding the
predictive accuracy of the algorithm as a function $\sqrt{\gamma/\varphi}.$
There it is clear that the characteristic time span where DF and AD compete is
not the day or the century, but clearly the year. To construct these figures
we have calculated the forward rate curves for different values of
$\sqrt{\gamma/\varphi}$ when one bond is removed from the constraining
dataset. The price of this missing bond is used afterwards as a benchmark to
test the accuracy of the resulting curves. Since we are interested in the
statistical performance we have done such comparison for 335 consecutive
trading days starting on Wednesday, November 08, 2000 and ending on Thursday,
March 07, 2002.

We are also interested in studying the impact of spreads in the constraining
dataset over the predictive accuracy. Therefore we present our results for
three spread patterns, namely constant spreads of 0\%, 0.5\% and 1\% in the
constraining dataset. Fig.(\ref{Graph4}) concentrates on the predictive
accuracy for the first 9 bonds of Table (\ref{Table2}) and Fig.(\ref{Graph5})
presents the same analysis for the remaining 2 bonds of Table (\ref{Table2}).
These last 2 bonds are the ones with the larger maturities in the complete
dataset. In particular for the last one with the largest maturity we have to
decide upon the methodology to extrapolate the forward curve outside the range
of the constraining dataset\footnote{By the range of the constraining dataset
we mean the range in time to maturity that goes from present to the last
maturity within the dataset.}. Therefore in Fig.(\ref{Graph5}) we present the
prediction accuracy using constant extrapolation from the last maturity in the
constraining dataset and \textquotedblleft$W$-generated\textquotedblright%
\ extrapolation that consists in utilizing the $W$-optimal forward rate curve
even outside the range of the constraining dataset. 

\section{Results and conclusions}

\label{sec conclusions} In this article we have presented a non-linear dynamic
programming algorithm designed for the calculation of positive definite
forward rate curves using data with or without spreads. We have included
multiple details of the algorithm aiming at practitioners not familiar with
the techniques of dynamic programming. We have illustrated the results of this
algorithms using the Swedish bond data for a one-parametric family of
smoothness measures.

The results and conclusions are the following:

\begin{itemize}
\item The proposed algorithm calculates forward rate curves in real time and
admits, without time or complexity penalizations, the use of any non-linear
\emph{local} objective function. Since it also handles non-linear constraints
it is possible to include within the constraining dataset any derivative
products with prices bearing some dependence on forward rates.

\item This is the first algorithm proposed in the literature that implements
the positivity constraint in the maximum smoothness framework. To the
knowledge of the authors the only other work that implements such constraint
outside this framework is the one of Wets et al. \cite{Wets}. The proposal in
\cite{Wets} has the advantage of using simple linear programming but do not
consider the presence of spreads minimizing instead the sum over the
\emph{modulus} of the difference between calculated prices and market prices.

\item For the objective functions and constraining datasets like the ones we
have used or more generally for the ones studied in \cite{Kwon} the algorithm
proposed in that reference offers better computing times at the expense of
ignoring the positivity constraint. Essentially the complexity in \cite{Kwon}
is $c^{3}$ and in ours is $n$ $c^{2}$ where $n$ is the number of time steps
and $c$ the number of constraints. For that reason, when this class of
objective functions and constraining datasets are used, a well coded algorithm
might try first the proposal in \cite{Kwon} (improving its treatment of the
spreads using e.g. log-barriers) and later, only if the result is not positive
definite, use our approach.

\item Since the optimization problem we are solving is non-convex we can not
discard the presence of several local minima (this is independent of the
presence of the positivity constraint \cite{Kwon}). However, we have tested
our algorithm starting from different seeds and in all cases we have arrived
at the same minima. These tests included hundreds of searches with initial
log-prices given by $\rho_{j}^{\left[  0\right]  }=\rho_{j}^{b}+x\left(
\rho_{j}^{a}-\rho_{j}^{b}\right)  $ (with $x\in\left[  0,1\right]  $ a flat
random variable) and initial forward rates given by several constant and
oscillating functions. With this comment we only want to convey our practical
experience and by no means we intend to say that we have exhaustively explored
the presence or absence of local minima.

\item It is clear that one way to avoid negative rates is just by increasing
spreads by hand. The advantage of using real spreads at a given moment is that
one can be sure of being consistent with market prices. If one observes in
Fig.(\ref{Graph3}) the large forward rate variations taking place for
different spread patterns no doubts should remain about the relevance of a
careful treatment of this issue.

\item From Fig.(\ref{Graph5}) we conclude that the inclusion of spreads can
remarkably improve the accuracy of resulting forward curves in the prediction
of market prices of long maturing bonds. In \cite{Kwon} it was pointed out
that the inclusion of spreads notably improved the smoothness of the forward
curve. To the knowledge of the authors this is the first time it is shown that
their presence also improves the prediction accuracy. For that reason we
believe that spreads should be considered even when the market data does not
provide such information. In that case the approach should consist in using a
cross validation technique to asses the optimal spread minimizing an error
criterion based in the prediction of market prices (for cross validation
methods see for example \cite{Weiss, Stone}). One possibility is using the
well know \textquotedblleft leave-one-out\textquotedblright\ cross-validation
to select the optimal spread much in the spirit suggested by
Figs.(\ref{Graph4}) and (\ref{Graph5}). Given a constraining set of $k$
products, the method actually consists in obtaining $k$ forward rate curves
for a given spread, each curve leaving out one of the constraining products
(bonds in our case) but using only the $k$ omitted products to compute an
error criterion like $%
{\displaystyle\sum_{i}^{k}}
\left\vert P_{i}^{pred.}-P_{i}^{market}\right\vert /P_{i}^{market}.$ Thus we
obtain a quantitative criterion to select an optimal level of
spread\footnote{This technique is stated here only to illustrate the existence
of an optimal non vanishing level of spread whenever no real spreads are
provided by the market. We do not intend to provide an efficient procedure to
find this level.}. For example for our dataset and our family of models
(measures) it can be conjectured from Figs.(\ref{Graph4}) and (\ref{Graph5})
that this optimal level of spread is typically around 0.5\%.

\item Figs.(\ref{Graph4}) and (\ref{Graph5}) strongly suggest that for low
spread patterns DF measure is more accurate than AD one. For bigger spreads
results do not clearly favor any particular measure.
\end{itemize}

\section*{Acknowledgements}

J. M. acknowledges the financial support from \emph{Tekniska H\"{o}gskolan,
Inst. f\"{o}r Fysik och M\"{a}tteknik}, Link\"{o}ping University, thanks J.
Shaw (Barclays) for calling his attention to ref. \cite{Wets}, F. Delbaen for
ref. \cite{Delbaen} and J. M. Eroles (HSBC) for proof reading the manuscript.
J. M. would also like to thank Dr. Per Olov Lindberg for the hospitality
during his stay at the Division of Optimization; Link\"{o}ping University. The
authors would like to thank the referees for their helpful comments on the manuscript.

\section*{Appendix A. Newton steps}

Since the objective function $Z$ in Eq.(\ref{total_obj}) is a non-quadratic
function of $f_{r}$, $\rho_{j}$ we will use an iterative quadratic
approximation (Newton steps) to find its minima. We start with feasible seeds
$f_{r}^{\left[  0\right]  }$ and $\rho_{j}^{\left[  0\right]  }$ fulfilling
inequality constraints (\ref{ineq_constr}) and we set up initial log-barriers
coefficients $\mu^{\left[  0\right]  }>0$ and $\tilde{\mu}^{\left[  0\right]
}>0.$ For newton iteration number $s$ we define
\begin{align}
\hat{f}_{r}^{\left[  s\right]  }  &  :=f_{r}^{\left[  s-1\right]  }+\Delta
_{r}^{\left[  s\right]  },\nonumber\\
\hat{\rho}_{j}^{\left[  s\right]  }  &  :=\rho_{j}^{\left[  s-1\right]
}+\sigma_{j}^{\left[  s\right]  }, \label{dynamic_sol}%
\end{align}
with $s=1,2,\cdots.$ The hat over $f$ and $\rho$ indicate that at each step
$s$, $\hat{f}^{\left[  s\right]  }$ and $\hat{\rho}^{\left[  s\right]  }$ are
the minima of the quadratic approximation and may not fulfill the inequality
constraints (\ref{ineq_constr}). To assure constraints (\ref{ineq_constr}) are
fulfilled a final redefinition $\hat{f}^{\left[  s\right]  }\rightarrow
f^{\left[  s\right]  },$ $\hat{\rho}^{\left[  s\right]  }\rightarrow
\rho^{\left[  s\right]  }$ is necessary after each Newton step. This
redefinition is explained in appendix C. Expanding $Z$ up to second order in
$\Delta_{r}^{\left[  s\right]  }$ and $\sigma_{j}^{\left[  s\right]  }$ we
write%
\begin{align}
Z^{\left[  s\right]  }  &  :=\left.  Z\left(  \hat{f}^{\left[  s\right]
},\hat{\rho}^{\left[  s\right]  },\lambda^{\left[  s\right]  }\right)
\right\vert _{O\left(  2\right)  }\nonumber\\
&  =\frac{1}{2}\Delta^{\left[  s\right]  T}Q^{\left[  s\right]  }%
\Delta^{\left[  s\right]  }+\Delta^{\left[  s\right]  T}B^{\left[  s\right]
}\lambda^{\left[  s\right]  }+\Delta^{\left[  s\right]  T}C^{\left[  s\right]
}+\lambda^{\left[  s\right]  T}a^{\left[  s\right]  }\nonumber\\
&  +\sigma^{\left[  s\right]  T}\lambda^{\left[  s\right]  }+\frac{1}{2}%
\sigma^{\left[  s\right]  T}M^{\left[  s\right]  }\sigma^{\left[  s\right]
}+\sigma^{\left[  s\right]  T}D^{\left[  s\right]  }+b^{\left[  s\right]  },
\label{quadratic}%
\end{align}
where $b^{\left[  s\right]  }$ collects all terms not depending on
$\Delta^{\left[  s\right]  }$, $\sigma^{\left[  s\right]  }$ or $\lambda
^{\left[  s\right]  }$. We will use square brackets around Newton step indices
and parenthesis around dynamic programming ones. Let us now work-out the
matrices involved in the quadratic approximation. From Eqs.(\ref{objective})
and (\ref{quadratic}) we immediately obtain%
\begin{align*}
&  \frac{1}{2}\Delta^{\left[  s\right]  T}Q^{\left[  s\right]  }%
\Delta^{\left[  s\right]  }+\Delta^{\left[  s\right]  T}C^{\left[  s\right]
}\\
&  =\gamma\sum_{r=1}^{n-1}\left(  \frac{f_{r+1}^{\left[  s-1\right]  }%
-f_{r}^{\left[  s-1\right]  }}{\xi_{r}}\right)  \left(  \frac{\Delta
_{r+1}^{\left[  s\right]  }-\Delta_{r}^{\left[  s\right]  }}{\xi_{r}}\right)
\xi_{r}+\frac{\gamma}{2}\sum_{r=1}^{n-1}\left(  \frac{\Delta_{r+1}^{\left[
s\right]  }-\Delta_{r}^{\left[  s\right]  }}{\xi_{r}}\right)  ^{2}\xi_{r}\\
&  +4\varphi\sum_{r=2}^{n-1}\left(  \frac{1}{\left(  \xi_{r-1}+\xi_{r}\right)
}\left(  \frac{1}{\xi_{r}}f_{r+1}^{\left[  s-1\right]  }+\frac{1}{\xi_{r-1}%
}f_{r-1}^{\left[  s-1\right]  }\right)  -\frac{1}{\xi_{r}\xi_{r-1}}%
f_{r}^{\left[  s-1\right]  }\right) \\
&  \times\left(  \frac{1}{\left(  \xi_{r-1}+\xi_{r}\right)  }\left(  \frac
{1}{\xi_{r}}\Delta_{r+1}^{\left[  s\right]  }+\frac{1}{\xi_{r-1}}\Delta
_{r-1}^{\left[  s\right]  }\right)  -\frac{1}{\xi_{r}\xi_{r-1}}\Delta
_{r}^{\left[  s\right]  }\right)  \xi_{r}\\
&  +2\varphi\sum_{r=2}^{n-1}\left(  \frac{1}{\left(  \xi_{r-1}+\xi_{r}\right)
}\left(  \frac{1}{\xi_{r}}\Delta_{r+1}^{\left[  s\right]  }+\frac{1}{\xi
_{r-1}}\Delta_{r-1}^{\left[  s\right]  }\right)  -\frac{1}{\xi_{r}\xi_{r-1}%
}\Delta_{r}^{\left[  s\right]  }\right)  ^{2}\xi_{r},
\end{align*}
and defining
\begin{align}
Y_{\lambda}^{\left[  s\right]  }  &  :=\sum_{j=1}^{m}\lambda_{j}^{\left[
s\right]  }\left[  \hat{\rho}_{j}^{\left[  s\right]  }-\ln\left(
v_{j}^{\left[  s\right]  }\right)  +\sum_{r=1}^{R_{n_{j}}^{\left(  j\right)
}-1}\hat{f}_{r}^{\left[  s\right]  }\xi_{r}\right]  ,\label{cosntr}\\
Y_{\mu}^{\left[  s\right]  }  &  :=-\mu^{\lbrack s-1]}\sum_{r=1}^{n}\ln\left(
\hat{f}_{r}^{\left[  s\right]  }\right)  ,\label{barrier1}\\
Y_{\tilde{\mu}}^{\left[  s\right]  }  &  :=-\tilde{\mu}^{\left[  s-1\right]
}\sum_{j=1}^{m}\left(  \ln\left(  \rho_{j}^{\left[  s\right]  }-\rho_{j}%
^{b}\right)  +\ln\left(  \rho_{j}^{a}-\rho_{j}^{\left[  s\right]  }\right)
\right)  , \label{barrier2}%
\end{align}
we have
\begin{align*}
Y_{\lambda}^{\left[  s\right]  }  &  =\sum_{j=1}^{m}\lambda_{j}^{\left[
s\right]  }\left[  \rho_{j}^{\left[  s-1\right]  }+\sigma_{j}^{\left[
s\right]  }-\ln\left(  v_{j}^{\left[  s-1\right]  }\right)  +\sum
_{r=1}^{R_{n_{j}}^{\left(  j\right)  }-1}\left(  f_{r}^{\left[  s-1\right]
}+\Delta_{r}^{\left[  s\right]  }\right)  \xi_{r}\right. \\
&  \left.  -\frac{1}{v_{j}^{\left[  s-1\right]  }}\sum\limits_{i=1}^{n_{j}%
-1}\sum_{r=R_{i}^{\left(  j\right)  }}^{R_{n_{j}}^{\left(  j\right)  }%
-1}\alpha_{ij}\exp\left(  \sum_{z=R_{i}^{\left(  j\right)  }}^{R_{n_{j}%
}^{\left(  j\right)  }-1}f_{z}^{\left[  s-1\right]  }\xi_{z}\right)
\Delta_{r}^{\left[  s\right]  }\xi_{r}+O\left(  \Delta_{r}^{\left[  s\right]
}{}^{2}\right)  \right]  ,\\
Y_{\mu}^{\left[  s\right]  }  &  =-\mu^{\left[  s-1\right]  }\sum_{r=1}%
^{n}\left(  \frac{\Delta_{r}^{\left[  s\right]  }}{f_{r}^{\left[  s-1\right]
}}-\frac{1}{2}\left(  \frac{\Delta_{r}^{\left[  s\right]  }}{f_{r}^{\left[
s-1\right]  }}\right)  ^{2}+O\left(  \Delta_{r}^{\left[  s\right]  }{}%
^{3}\right)  \right)  ,\\
Y_{\tilde{\mu}}^{\left[  s\right]  }  &  =-\tilde{\mu}^{\left[  s-1\right]
}\sum_{j=1}^{m}\left(  \ln\left(  \rho_{j}^{\left[  s-1\right]  }-\rho_{j}%
^{b}\right)  +\ln\left(  \rho_{j}^{a}-\rho_{j}^{\left[  s-1\right]  }\right)
\right) \\
&  +\tilde{\mu}^{\left[  s-1\right]  }\sum_{j=1}^{m}\left(  \frac{1}{\rho
_{j}^{a}-\rho_{j}^{\left[  s-1\right]  }}-\frac{1}{\rho_{j}^{\left[
s-1\right]  }-\rho_{j}^{b}}\right)  \sigma_{j}^{\left[  s\right]  }\\
&  +\frac{\tilde{\mu}^{\left[  s-1\right]  }}{2}\sum_{j=1}^{m}\left[  \frac
{1}{\left(  \rho_{j}^{\left[  s-1\right]  }-\rho_{j}^{a}\right)  ^{2}}%
+\frac{1}{\left(  \rho_{j}^{\left[  s-1\right]  }-\rho_{j}^{b}\right)  ^{2}%
}\right]  \sigma_{j}^{\left[  s\right]  2}+O\left(  \sigma_{j}^{\left[
s\right]  3}\right)  .
\end{align*}
Thus defining
\[
\chi\left(  r,x,y\right)  :=\left\{
\begin{array}
[c]{cc}%
1 & x\leq r\leq y\\
0 & \mathrm{otherwise}%
\end{array}
\right.  ,\qquad\delta_{i,j}:=\left\{
\begin{array}
[c]{cc}%
1 & i=j\\
0 & \mathrm{otherwise}%
\end{array}
\right.  ,
\]
from above expansions and Eq.(\ref{quadratic}) we immediately obtain%
\begin{align*}
B_{r,j}^{\left[  s\right]  }  &  =\xi_{r}\chi\left(  r,1,R_{n_{j}}^{\left(
j\right)  }-1\right)  -\frac{1}{v_{j}^{\left[  s-1\right]  }}\sum
\limits_{i=1}^{n_{j}-1}\alpha_{ij}\exp\left(  \sum_{z=R_{i}^{\left(  j\right)
}}^{R_{n_{j}}^{\left(  j\right)  }-1}f_{z}^{\left[  s-1\right]  }\xi
_{z}\right)  \xi_{r}\chi\left(  r,R_{i}^{\left(  j\right)  },R_{n_{j}%
}^{\left(  j\right)  }-1\right)  ,\\
a_{j}^{\left[  s\right]  }  &  =\rho_{j}^{\left[  s-1\right]  }-\ln\left(
v_{j}^{\left[  s-1\right]  }\right)  +\sum_{r=1}^{R_{n_{j}}^{\left(  j\right)
}-1}f_{r}^{\left[  s-1\right]  }\xi_{r},\\
Q_{r,x}^{\left[  s\right]  }  &  =4\varphi\left[  \left(  \frac{\left(
1-\delta_{r,1}\right)  \left(  1-\delta_{r,2}\right)  \xi_{r-1}}{\left(
\xi_{r-2}+\xi_{r-1}\right)  ^{2}\xi_{r-1}^{2}}+\frac{\left(  1-\delta
_{r,n}\right)  \left(  1-\delta_{r,1}\right)  \xi_{r}}{\xi_{r}^{2}\xi
_{r-1}^{2}}+\frac{\left(  1-\delta_{r,n-1}\right)  \left(  1-\delta
_{r,n}\right)  \xi_{r+1}}{\left(  \xi_{r}+\xi_{r+1}\right)  ^{2}\xi_{r}^{2}%
}\right)  \delta_{r,x}\right. \\
&  +\frac{\xi_{r+1}\delta_{r+2,x}}{\left(  \xi_{r}+\xi_{r+1}\right)  ^{2}%
\xi_{r+1}\xi_{r}}+\frac{\xi_{x+1}\delta_{r-2,x}}{\left(  \xi_{x}+\xi
_{x+1}\right)  ^{2}\xi_{x+1}\xi_{x}}-\frac{\left(  1-\delta_{r,1}\right)
\xi_{r}\delta_{r+1,x}}{\left(  \xi_{r-1}+\xi_{r}\right)  \xi_{r-1}\xi_{r}^{2}%
}\\
&  \left.  -\frac{\left(  1-\delta_{x,1}\right)  \xi_{x}\delta_{r-1,x}%
}{\left(  \xi_{x-1}+\xi_{x}\right)  \xi_{x-1}\xi_{x}^{2}}-\frac{\left(
1-\delta_{r,n}\right)  \xi_{r}\delta_{r-1,x}}{\left(  \xi_{r-1}+\xi
_{r}\right)  \xi_{r}\xi_{r-1}^{2}}-\frac{\left(  1-\delta_{x,n}\right)
\xi_{x}\delta_{r+1,x}}{\left(  \xi_{x-1}+\xi_{x}\right)  \xi_{x}\xi_{x-1}^{2}%
}\right] \\
&  +\gamma\left[  \left(  1-\delta_{r,n}\right)  \frac{\xi_{r}}{\xi_{r}^{2}%
}\delta_{r,x}+\left(  1-\delta_{r,1}\right)  \frac{\xi_{r-1}}{\xi_{r-1}^{2}%
}\delta_{r,x}-\frac{\xi_{r}}{\xi_{r}^{2}}\delta_{r+1,x}-\frac{\xi_{x}}{\xi
_{x}^{2}}\delta_{r-1,x}\right]  +\frac{\mu^{\left[  s-1\right]  }}%
{f_{r}^{\left[  s-1\right]  2}}\delta_{r,x},\\
C_{r}^{\left[  s\right]  }  &  =-\frac{\mu^{\left[  s-1\right]  }}%
{f_{r}^{\left[  s-1\right]  }}+\gamma\frac{f_{r}^{\left[  s-1\right]
}-f_{r-1}^{\left[  s-1\right]  }}{\xi_{r-1}}\frac{\left(  1-\delta
_{r,1}\right)  }{\xi_{r-1}}\xi_{r-1}-\gamma\frac{f_{r+1}^{\left[  s-1\right]
}-f_{r}^{\left[  s-1\right]  }}{\xi_{r}}\frac{\left(  1-\delta_{r,n}\right)
}{\xi_{r}}\xi_{r}\\
&  +4\varphi\left(  \frac{f_{r}^{\left[  s-1\right]  }-f_{r-1}^{\left[
s-1\right]  }}{\xi_{r-1}}+\frac{f_{r-2}^{\left[  s-1\right]  }-f_{r-1}%
^{\left[  s-1\right]  }}{\xi_{r-2}}\right)  \frac{\left(  1-\delta
_{r,1}\right)  \left(  1-\delta_{r,2}\right)  }{\left(  \xi_{r-2}+\xi
_{r-1}\right)  ^{2}\xi_{r-1}}\xi_{r-1}\\
&  +4\varphi\left(  \frac{f_{r+2}^{\left[  s-1\right]  }-f_{r+1}^{\left[
s-1\right]  }}{\xi_{r+1}}+\frac{f_{r}^{\left[  s-1\right]  }-f_{r+1}^{\left[
s-1\right]  }}{\xi_{r}}\right)  \frac{\left(  1-\delta_{r,n}\right)  \left(
1-\delta_{r,n-1}\right)  }{\left(  \xi_{r}+\xi_{r+1}\right)  ^{2}\xi_{r}}%
\xi_{r+1}\\
&  -4\varphi\left(  \frac{f_{r+1}^{\left[  s-1\right]  }-f_{r}^{\left[
s-1\right]  }}{\xi_{r}}+\frac{f_{r-1}^{\left[  s-1\right]  }-f_{r}^{\left[
s-1\right]  }}{\xi_{r-1}}\right)  \frac{\left(  1-\delta_{r,1}\right)  \left(
1-\delta_{r,n}\right)  }{\left(  \xi_{r-1}+\xi_{r}\right)  \xi_{r}\xi_{r-1}%
}\xi_{r},\\
M_{j,k}^{\left[  s\right]  }  &  =\tilde{\mu}^{\left[  s-1\right]  }\left[
\frac{1}{\left(  \rho_{j}^{a}-\rho_{j}^{\left[  s-1\right]  }\right)  ^{2}%
}+\frac{1}{\left(  \rho_{j}^{b}-\rho_{j}^{\left[  s-1\right]  }\right)  ^{2}%
}\right]  \delta_{j,k},\\
D_{j}^{\left[  s\right]  }  &  =\tilde{\mu}^{\left[  s-1\right]  }\left(
\frac{1}{\rho_{j}^{a}-\rho_{j}^{\left[  s-1\right]  }}+\frac{1}{\rho_{j}%
^{b}-\rho_{j}^{\left[  s-1\right]  }}\right)  .
\end{align*}

\section*{Appendix B. Dynamic Programming for a quadratic objective}

Given positive definite symmetric matrices $Q\in R^{n,n}$ and $M\in R^{c,c}$
and the $n$-vector $C$ and $c$-vector $D$, we want to minimize the objective
function
\[
W\left(  \Delta,\sigma\right)  :=\frac{1}{2}\Delta^{T}Q\Delta+\Delta
^{T}C+\frac{1}{2}\sigma^{T}M\sigma+\sigma^{T}D,
\]
subject to the set of constraints%
\begin{equation}
a_{j}+\sigma_{j}+\sum_{k=0}^{n}\Delta_{k}B_{k,j}=0,\qquad j=1,\ldots,c,
\label{constraints}%
\end{equation}
thus defining a new objective function
\[
Z\left(  \Delta,\sigma,\lambda\right)  :=W\left(  \Delta,\sigma\right)
+\left(  \Delta^{T}B+a^{T}+\sigma^{T}\right)  \lambda.
\]
For any $Q$ we define%
\[
d:=\max\limits_{k}\left(  \sum\limits_{s=1}^{k-1}\delta\left(  Q_{k,s}\right)
\right)  ,\qquad\delta\left(  x\right)  :=\left\{
\begin{array}
[c]{cc}%
1 & x\neq0\\
0 & x=0
\end{array}
\right.  ,
\]
where e.g. $d=0$, $1$ for diagonal, tridiagonal matrices respectively. When
$d<<n$ and $c<<n$ then $Z$ can be solved efficiently with dynamic programming.
To set up the notation and for those readers not familiar with dynamic
programming let us here briefly explain the basics of this well known method
\cite{dynamic}. We start defining
\begin{align}
Z^{\left(  n\right)  }  &  :=Z,\nonumber\\
Z^{\left(  q-1\right)  }  &  :=\left.  Z^{\left(  q\right)  }\right\vert
_{\Delta_{q}=\Delta_{q}^{\ast}},\qquad q=1,\ldots,n, \label{fundamental}%
\end{align}
where $\Delta_{q}^{\ast}$ satisfies%
\begin{equation}
\left.  \triangledown_{\Delta_{q}}Z^{\left(  q\right)  }\right\vert
_{\Delta_{q}=\Delta_{q}^{\ast}}=0,\qquad q=1,\ldots,n. \label{extrema}%
\end{equation}
For $q=0,\ldots,n$ we use the inductive hypothesis%
\begin{align}
Z^{\left(  q\right)  }  &  =\frac{1}{2}\Delta^{T}Q^{\left(  q\right)  }%
\Delta+\Delta^{T}B^{\left(  q\right)  }\lambda+\Delta^{T}C^{\left(  q\right)
}+\lambda^{T}a^{\left(  q\right)  }\nonumber\\
&  +\frac{1}{2}\lambda^{T}G^{\left(  q\right)  }\lambda+b^{\left(  q\right)
}+\sigma^{T}\lambda+\frac{1}{2}\sigma^{T}M\sigma+\sigma^{T}D,
\label{hypothesis}%
\end{align}
where $Q^{\left(  q\right)  }\in R^{q,q},$ $B^{\left(  q\right)  }\in R^{c,q}$
and $G^{\left(  q\right)  }\in R^{c,c}$. The final step of this backwards
process consists in obtaining $\sigma^{\ast}$ and $\lambda^{\ast}$ satisfying%
\begin{align}
\left.  \triangledown_{\sigma}Z^{\left(  0\right)  }\right\vert _{\sigma
=\sigma^{\ast}}  &  =0,\nonumber\\
\left.  \triangledown_{\lambda}Z^{\left(  0\right)  }\right\vert
_{\lambda=\lambda^{\ast}}  &  =0. \label{Lagrange}%
\end{align}
From Eqs.(\ref{extrema}) and (\ref{hypothesis}) we obtain
\begin{equation}
\Delta_{q}^{\ast}=-\frac{1}{Q_{q,q}^{\left(  q\right)  }}\left(
\sum\limits_{r=q-d}^{q-1}Q_{q,r}^{\left(  q\right)  }\Delta_{r}+B_{q}^{\left(
q\right)  }\lambda+C_{q}^{\left(  q\right)  }\right)  , \label{first}%
\end{equation}
and plugging Eq.(\ref{first}) into Eq.(\ref{fundamental}) we obtain%
\begin{align*}
Z^{\left(  q-1\right)  }  &  =\left.  Z^{\left(  q\right)  }\right\vert
_{\Delta_{q}=0}-\frac{1}{2}\frac{1}{Q_{q,q}^{\left(  q\right)  }}\left(
\sum\limits_{r=q-d}^{q-1}Q_{q,r}^{\left(  q\right)  }\Delta_{r}+B_{q}^{\left(
q\right)  }\lambda+C_{q}^{\left(  q\right)  }\right)  ^{2}\\
&  =\frac{1}{2}\sum\limits_{\left\{  r,s\right\}  =1}^{q-1}\Delta_{r}\left(
Q_{r,s}^{\left(  q\right)  }-\frac{Q_{q,r}^{\left(  q\right)  }Q_{q,s}%
^{\left(  q\right)  }}{Q_{q,q}^{\left(  q\right)  }}\theta\left(
r-q+d\right)  \theta\left(  s-q+d\right)  \right)  \Delta_{s}\\
&  +\sum\limits_{r=1}^{q-1}\Delta_{r}\left(  \left(  B_{r}^{\left(  q\right)
}-\frac{Q_{q,r}^{\left(  q\right)  }B_{q}^{\left(  q\right)  }}{Q_{q,q}%
^{\left(  q\right)  }}\theta\left(  r-q+d\right)  \right)  \lambda
+C_{r}^{\left(  q\right)  }-\frac{Q_{q,r}^{\left(  q\right)  }C_{q}^{\left(
q\right)  }}{Q_{q,q}^{\left(  q\right)  }}\theta\left(  r-q+d\right)  \right)
\\
&  +\frac{1}{2}\lambda^{T}\left(  G^{\left(  q\right)  }-\frac{B_{q}^{\left(
q\right)  T}B_{q}^{\left(  q\right)  }}{Q_{q,q}^{\left(  q\right)  }}\right)
\lambda+\left(  \alpha^{\left(  q\right)  T}-\frac{C_{q}^{\left(  q\right)
}B_{q}^{\left(  q\right)  }}{Q_{q,q}^{\left(  q\right)  }}\right)
\lambda-\frac{1}{2}\frac{C_{q}^{\left(  q\right)  2}}{Q_{q,q}^{\left(
q\right)  }}+b^{\left(  q\right)  }\\
&  +\sigma^{T}\lambda+\frac{1}{2}\sigma^{T}M\sigma+\sigma^{T}D,
\end{align*}
hence obtaining%
\begin{align}
Q_{r,s}^{\left(  q-1\right)  }  &  =Q_{r,s}^{\left(  q\right)  }-\frac
{Q_{q,r}^{\left(  q\right)  }Q_{q,s}^{\left(  q\right)  }}{Q_{q,q}^{\left(
q\right)  }}\theta\left(  r-q+d\right)  \theta\left(  s-q+d\right)
,\label{Q}\\
B_{r}^{\left(  q-1\right)  }  &  =B_{r}^{\left(  q\right)  }-\frac
{Q_{q,r}^{\left(  q\right)  }B_{q}^{\left(  q\right)  }}{\tilde{Q}%
_{q,q}^{\left(  q\right)  }}\theta\left(  r-q+d\right)  ,\label{B}\\
C_{r}^{\left(  q-1\right)  }  &  =C_{r}^{\left(  q\right)  }-\frac
{Q_{q,r}^{\left(  q\right)  }C_{q}^{\left(  q\right)  }}{Q_{q,q}^{\left(
q\right)  }}\theta\left(  r-q+d\right)  ,\label{C}\\
G^{\left(  q-1\right)  }  &  =G^{\left(  q\right)  }-\frac{B_{q}^{\left(
q\right)  T}B_{q}^{\left(  q\right)  }}{Q_{q,q}^{\left(  q\right)  }%
},\label{G}\\
a^{\left(  q-1\right)  }  &  =a^{\left(  q\right)  }-\frac{C_{q}^{\left(
q\right)  }B_{q}^{\left(  q\right)  }}{Q_{q,q}^{\left(  q\right)  }}%
,\label{a}\\
b^{\left(  q-1\right)  }  &  =b^{\left(  q\right)  }-\frac{1}{2}\frac
{C_{q}^{\left(  q\right)  2}}{Q_{q,q}^{\left(  q\right)  }}, \label{b}%
\end{align}
where%
\[
\theta\left(  x\right)  :=\left\{
\begin{array}
[c]{cc}%
1 & x\geq0,\\
0 & x<0.
\end{array}
\right.
\]

Now from Eq.(\ref{Lagrange}) we obtain%

\begin{align*}
M\sigma^{\ast}+\lambda^{\ast}+D  &  =0,\\
G^{\left(  0\right)  }\lambda^{\ast}+\sigma^{\ast}+a^{\left(  0\right)  }  &
=0,
\end{align*}
or%
\begin{align}
\lambda^{\ast}  &  =\left(  M^{-1}-G^{\left(  0\right)  }\right)  ^{-1}\left(
a^{\left(  0\right)  }-M^{-1}D\right)  =\left(  I-MG^{\left(  0\right)
}\right)  ^{-1}\left(  Ma^{\left(  0\right)  }-D\right)  ,\nonumber\\
\sigma^{\ast}  &  =-M^{-1}\left(  \lambda^{\ast}+D\right)  , \label{sig}%
\end{align}
and using Eq.(\ref{G}) we obtain%
\[
\lambda^{\ast}=\left(  M^{-1}+\sum\limits_{q=1}^{n}\frac{B_{q}^{\left(
q\right)  T}B_{q}^{\left(  q\right)  }}{Q_{q,q}^{\left(  q\right)  }}\right)
^{-1}\left(  a^{\left(  0\right)  }-M^{-1}D\right)  =\left(  I+M\sum
\limits_{q=1}^{n}\frac{B_{q}^{\left(  q\right)  T}B_{q}^{\left(  q\right)  }%
}{Q_{q,q}^{\left(  q\right)  }}\right)  ^{-1}\left(  Ma^{\left(  0\right)
}-D\right)  .
\]
Finally from $\lambda^{\ast}$ and Eq.(\ref{sig}) we obtain $\sigma^{\ast} $
and then using Eqs.(\ref{first}) and (\ref{Q}-\ref{C}) we obtain $\Delta
^{\ast}$ moving forward in the $q$ index.

\subsection*{B.1 Computing times of a dynamic programming iteration}

The time necessary to compute all $Q^{(i)}$ is proportional to%
\[
T_{1}=\left(  n-d\right)  \frac{d\left(  d+1\right)  }{2}+\sum\limits_{i=1}%
^{d-1}\frac{i\left(  i+1\right)  }{2}=\frac{d\left(  d+1\right)  }{2}\left(
n-\frac{2d+1}{3}\right)  .
\]
To calculate the inverse of $M^{-1}-G^{\left(  0\right)  }$, that is a
$c\times c$ symmetric matrix, we require a computing time proportional to%
\[
T_{2}=\frac{1}{6}\left(  c^{3}-c\right)  .
\]
All $C^{(i)},B^{(i)}$ and $a^{(i)}$ require, respectively, computing times
proportional to%
\begin{align*}
T_{3}  &  =\left(  n-d\right)  d+\sum\limits_{i=1}^{d-1}i=\frac{d}{2}\left(
2n-d-1\right)  ,\\
T_{4}  &  =cT_{3},\\
T_{5}  &  =cn.
\end{align*}
Finally the calculation of $G^{\left(  0\right)  }$ requires a computing time
proportional to
\[
T_{6}=\frac{c\left(  c+1\right)  }{2}n.
\]
In the forward rate calculation typically we have $n>c>d$ (in
Eq.(\ref{objective}) we have $d=2$) and therefore the maximum delay would be
given by $T_{6}.$

\section*{Appendix C. Fulfilling inequalities and updating log-barriers}

The Newton step obtained from Eq.(\ref{dynamic_sol}) may not fulfill the
inequality constraints (\ref{ineq_constr}). To satisfy such constraints we
search for a $\alpha^{\lbrack s]}$ in the interval $\left[  0,1\right]  $ such
that%
\begin{equation}
f_{r}^{\left[  s-1\right]  }+\alpha^{\lbrack s]}\Delta_{r}^{\left[  s\right]
}\geq0,\qquad\rho_{j}^{b}\leq\rho_{j}^{\left[  s-1\right]  }+\alpha^{\lbrack
s]}\sigma_{j}^{\left[  s\right]  }\leq\rho_{j}^{a}.\label{constr_alpha}%
\end{equation}
In order to do this we first determine the maximum $\alpha_{\max}^{[s]}$ in
the interval $\left[  0,1\right]  $ satisfying Eq.(\ref{constr_alpha}). That
is, given the sets of points%
\begin{align*}
A &  :=\left\{  r,\text{ }\hat{f}_{r}^{\left[  s\right]  }\leq0\right\}  ,\\
\tilde{A}_{\leq} &  :=\left\{  j,\text{ }\hat{\rho}_{j}^{\left[  s\right]
}\leq\rho_{j}^{b}\right\}  ,\\
\tilde{A}_{\geq} &  :=\left\{  j,\text{ }\hat{\rho}_{j}^{\left[  s\right]
}\geq\rho_{j}^{a}\right\}  ,
\end{align*}
we have%
\[
\alpha_{\max}^{[s]}=\min\left(  \min\limits_{r\in A}\left(  \frac
{-f_{r}^{\left[  s-1\right]  }}{\Delta_{r}^{\left[  s\right]  }}\right)
,\min\limits_{j\in\tilde{A}_{\leq}}\left(  \frac{\rho_{j}^{b}-\rho
_{j}^{\left[  s-1\right]  }}{\sigma_{j}^{\left[  s\right]  }}\right)
,\min\limits_{j\in\tilde{A}_{\geq}}\left(  \frac{\rho_{j}^{a}-\rho
_{j}^{\left[  s-1\right]  }}{\sigma_{j}^{\left[  s\right]  }}\right)  \right)
,
\]
and then we take%
\[
\alpha^{\lbrack s]}:=\beta\alpha_{\max}^{[s]},
\]
with $0<\beta<1$ (in our implementation we have taken $\beta=.9$ that is a
standard election in the optimization literature). Once we have $\alpha
^{\lbrack s]}$ we define%
\begin{align}
f_{r}^{\left[  s\right]  } &  =f_{r}^{\left[  s-1\right]  }+\alpha^{\lbrack
s]}\Delta_{r}^{\left[  s\right]  },\nonumber\\
\rho_{j}^{\left[  s\right]  } &  =\rho_{j}^{\left[  s-1\right]  }%
+\alpha^{\lbrack s]}\sigma_{j}^{\left[  s\right]  },\nonumber\\
\mu^{\left[  s\right]  } &  =\max(\Psi\left(  \alpha^{\lbrack s]}\right)
\mu^{\left[  s-1\right]  },\mu_{\min}),\nonumber\\
\tilde{\mu}^{\left[  s\right]  } &  =\max(\Psi\left(  \alpha^{\lbrack
s]}\right)  \tilde{\mu}^{\left[  s-1\right]  },\tilde{\mu}_{\min
}),\label{update}%
\end{align}
where $\mu_{\min}$ and $\tilde{\mu}_{\min}$ are positive small values
guaranteeing that matrices $Q$ and $M$ are positive definite and
$\Psi:[0,1]\rightarrow R^{+}$ is a monotonically decreasing function
satisfying $\Psi\left(  0\right)  =1$. In our implementation we have taken
$\Psi$ of the form%
\begin{equation}
\Psi\left(  \alpha\right)  =\left(  1-l\right)  \left(  1-\alpha\right)
^{\xi}+l,\label{factorf}%
\end{equation}
with $l=10^{-2}$, $\xi=1$. The value of $l$ controls how fast barriers are
reduced. In our implementation we have found that adequate values for $l$
range $10^{-3}\lesssim l\lesssim10^{-1}$. The value of $\xi$ controls the
non-linearity of $\Psi$. We have found that the simple linear response
provides good performance albeit convergence time is not significantly
affected for $\xi$ in the range $0.6\lesssim\xi\lesssim2$

In this way we iterate the algorithm $n_{it}$ times until a given termination
criterion is met. Defining%
\begin{align*}
\delta_{\ln W}^{\left[  s\right]  }  &  :=\ln\left(  W\left(  f^{[s]}\right)
\right)  -\ln\left(  W\left(  f^{[s-1]}\right)  \right)  ,\\
\epsilon^{\left[  s\right]  }  &  :=\max\limits_{j}\left(  \left\vert \rho
_{j}^{\left[  s\right]  }-\ln\left(  v_{j}^{\left[  s\right]  }\right)
+\sum_{r=1}^{R_{n_{j}}^{\left(  j\right)  }-1}f_{r}^{\left[  s\right]  }%
\xi_{r}\right\vert \right)  ,
\end{align*}
we have chosen the following termination criterion%
\begin{align}
&  ~n_{it}=N_{\max}~~\mathrm{or}~~\left\{  \left[  W^{[s]}<W_{zero}%
~~\mathrm{or}~~\left(  \left\vert \delta_{\ln W}^{\left[  s\right]
}\right\vert <\delta_{\ln W}^{\max}~~\mathrm{and}~~\left\vert \delta_{\ln
W}^{\left[  s-1\right]  }\right\vert <\delta_{\ln W}^{\max}\right)  \right]
\right. \nonumber\\
&  \left.  \mathrm{and}~~n_{it}>N_{\min}~~\mathrm{and}~~\mu^{\left[  s\right]
}<\mu_{\max}~~\mathrm{and}~~\tilde{\mu}^{\left[  s\right]  }<\tilde{\mu}%
_{\max}~~\mathrm{and}~~\epsilon^{\left[  s\right]  }<\epsilon_{\max}\right\}
. \label{criterion}%
\end{align}
In our implementation we have taken%
\begin{align}
\mu^{\lbrack0]}  &  =10^{-1},\quad\mu_{\min}=10^{-10},\quad\mu_{\max}%
=10^{-6},\nonumber\\
\tilde{\mu}^{[0]}  &  =10^{+1},\quad\tilde{\mu}_{\min}=10^{-10},\quad
\tilde{\mu}_{\max}=10^{-6},\nonumber\\
\epsilon_{\max}  &  =10^{-8},\quad\delta_{\ln W}^{\max}=10^{-2},\quad
W_{zero}=10^{-9},\nonumber\\
N_{\min}  &  =5,\quad N_{\max}=60. \label{values}%
\end{align}

Obviously there is considerable latitude to change the heuristic values
assigned to the above parameters. Let us finish this appendix making some
comments regarding their robustness.

$N_{\min}$ is there to guarantee a minimum number of iterations so as to have
$\delta_{\ln W}^{\left[  s\right]  }$ and $\delta_{\ln W}^{\left[  s-1\right]
}$ well defined and to avoid premature termination in the improbable case
where the other criteria incorrectly suggest convergence. For this purpose is
enough to take $2\leq N_{\min}\lesssim5$.

$N_{\max}$ serves as a maximum limit to secure termination even if convergence
is not achieved and therefore an alarm should be provided whenever
$n_{it}=N_{\max}$. From our experience we observe that is more than enough to
take $50\lesssim N_{\max}\lesssim100.$

$W_{zero}$ sets our precision to consider a given forward rate curve as a
straight line. The order of magnitude of $W_{zero}$ should be taken much lower
than the typical order of magnitude of the observed optimal $W.$ The value of
the optimal $W$ depends not only on the constraining data but also on $\gamma$
and $\varphi.$ We have adopted the practise of spanning the range $0\leq
\sqrt{\gamma/\varphi}\leq1$ year$^{-1}$ taking $0\leq\gamma\leq1$ year$^{3}$,
$\varphi=1$ year$^{5}$ and the range $1$ year$^{-1}\leq\sqrt{\gamma/\varphi
}\leq\infty$ taking $\gamma=1$ year$^{3}$, $0\leq\varphi\leq1$ year$^{5}$.
With this convention we have found reasonable to take $10^{-10}\lesssim
W_{zero}\lesssim10^{-8}$ for the full range $0\leq\sqrt{\gamma/\varphi}%
\leq\infty.$

$\delta_{\ln W}^{\max}$ sets the maximum variation of $\ln\left(  W\right)  $
between newton steps that is accepted before termination. Note that in
(\ref{criterion}) to have ($\left\vert \delta_{\ln W}^{\left[  s\right]
}\right\vert <\delta_{\ln W}^{\max}~~\mathrm{and}~~\left\vert \delta_{\ln
W}^{\left[  s-1\right]  }\right\vert <\delta_{\ln W}^{\max}$) is a necessary
but not sufficient condition for termination. Therefore sending $\delta_{\ln
W}^{\max}\rightarrow\infty$ has no major impact and is equivalent to rely only
on the error $\epsilon^{\left[  s\right]  }$ and the barrier coefficients as
indicators of convergence. On the contrary excessively reducing $\delta_{\ln
W}^{\max}$ can generate unnecessary iterations. We have tested that for values
satisfying $\delta_{\ln W}^{\max}\gtrsim10^{-5}$ we do not have any drastic
increase in convergence time.

$\epsilon_{\max}$ controls the error in the constraints and is the most
important parameter in (\ref{criterion}). A too large value of $\epsilon
_{\max}$ reduces the accuracy of the result and a too small value can give
rise to unnecessary iterations. We have observed acceptable results for
$\epsilon_{\max}$ in the range $10^{-4}\lesssim\epsilon_{\max}\lesssim
10^{-10}.$

$\mu_{\max}$ and $\tilde{\mu}_{\max}$ control the maximum allowed values for
the logarithmic barriers implementing the positivity and spread constraints
respectively. Large values for these parameters can make log barriers to have
a residual influence in the feasible region. Acceptable values of these
maximum weights range in $\mu_{\min}<\mu_{\max}\lesssim10^{-6}$ and
$\tilde{\mu}_{\min}<\tilde{\mu}_{\max}\lesssim10^{-6}.$ As explained above
keeping parameters $\mu_{\min}$ and $\tilde{\mu}_{\min}$ positive guarantees
that matrices $Q$ and $M$ are positive definite which in turn is sufficient to
guarantee the existence of an optimal solution in each iteration. Hence
$\mu_{\min}$ and $\tilde{\mu}_{\min}$ should be chosen as small as possible
without interfering with the numerical stability of the algorithm. Using
double precision in our program we have found the values given in
(\ref{criterion}) as a good compromise.

Finally $\mu^{\lbrack0]}$ and $\tilde{\mu}^{[0]}$ are the initial values for
the log barriers coefficients. Taking very large values for $\mu^{\lbrack0]}$
and $\tilde{\mu}^{[0]}$ increases the convergence time because we need more
time to reduce the barriers. Taking too small values for $\mu^{\lbrack0]}$ and
$\tilde{\mu}^{[0]}$ also increases the convergence time because time is wasted
exploring unfeasible solutions. Moreover, we have observed that convergence
and stability is improved if the contributions to the objective of the two
barrier terms are kept balanced. This is achieved setting $\tilde{\mu}%
^{[0]}\simeq\frac{n}{2m}\mu^{\lbrack0]}$ (see Eqs.(\ref{barrier1}) and
(\ref{barrier2})) and using the same updating factor $\Psi\left(
\alpha^{\lbrack s]}\right)  $ for $\mu^{\lbrack s]}$ and $\tilde{\mu}^{\left[
s\right]  }$ (see Eq.(\ref{update})). Keeping $\tilde{\mu}^{[0]}\simeq\frac
{n}{2m}\mu^{\lbrack0]}$ we have found that the time of convergence is stable
for $\mu^{\lbrack0]}$ in the range $10^{-4}\lesssim\mu^{\lbrack0]}%
\lesssim10^{+1}.$

\section*{Appendix D. Bond tables}

In this work we have used the following data tables and conventions.
\begin{table}[h]
\centering
\begin{tabular}
[c]{|l||l|l|l|l|l|l|l|l|l|l|l|}\hline
{\tiny Date}$\backslash${\tiny Bond (SO)} & {\tiny 1033} & {\tiny 1042} &
{\tiny 1035} & {\tiny 1044} & {\tiny 1038} & {\tiny 1037} & {\tiny 1040} &
{\tiny 1043} & {\tiny 1034} & {\tiny 1045} & {\tiny 1041}\\\hline\hline
{\tiny Friday 06} & {\tiny 4.86} & {\tiny 4.92} & {\tiny 5.06} & {\tiny 5.15}
& {\tiny 5.26} & {\tiny 5.27} & {\tiny 5.355} & {\tiny 5.4} & {\tiny 5.395} &
{\tiny 5.46} & {\tiny 5.655}\\\hline
{\tiny Monday 09} & {\tiny 4.905} & {\tiny 4.965} & {\tiny 5.11} & {\tiny 5.2}
& {\tiny 5.05} & {\tiny 5.325} & {\tiny 5.4} & {\tiny 5.455} & {\tiny 5.23} &
{\tiny 5.51} & {\tiny 5.7}\\\hline
{\tiny Tuesday 10} & {\tiny 4.885} & {\tiny 4.945} & {\tiny 5.095} &
{\tiny 5.185} & {\tiny 4.93} & {\tiny 5.295} & {\tiny 5.395} & {\tiny 5.435} &
{\tiny 5.24} & {\tiny 5.495} & {\tiny 5.68}\\\hline
{\tiny Wednesday 11} & {\tiny 4.865} & {\tiny 4.92} & {\tiny 5.06} &
{\tiny 5.15} & {\tiny 4.93} & {\tiny 5.275} & {\tiny 5.355} & {\tiny 5.41} &
{\tiny 5.415} & {\tiny 5.465} & {\tiny 5.65}\\\hline
{\tiny Thursday 12} & {\tiny 4.835} & {\tiny 4.885} & {\tiny 5.025} &
{\tiny 5.125} & {\tiny 4.93} & {\tiny 5.25} & {\tiny 5.355} & {\tiny 5.405} &
{\tiny 5.405} & {\tiny 5.465} & {\tiny 5.66}\\\hline
{\tiny Friday 13} & {\tiny 4.84} & {\tiny 4.89} & {\tiny 5.045} &
{\tiny 5.145} & {\tiny 4.93} & {\tiny 5.27} & {\tiny 5.38} & {\tiny 5.43} &
{\tiny 5.43} & {\tiny 5.495} & {\tiny 5.69}\\\hline
{\tiny Monday 16} & {\tiny 4.825} & {\tiny 4.87} & {\tiny 5.035} &
{\tiny 5.13} & {\tiny 4.93} & {\tiny 5.25} & {\tiny 5.36} & {\tiny 5.415} &
{\tiny 5.41} & {\tiny 5.475} & {\tiny 5.66}\\\hline
{\tiny Tuesday 17} & {\tiny 4.805} & {\tiny 4.85} & {\tiny 5.015} &
{\tiny 5.11} & {\tiny 4.93} & {\tiny 5.23} & {\tiny 5.34} & {\tiny 5.395} &
{\tiny 5.395} & {\tiny 5.455} & {\tiny 5.64}\\\hline
{\tiny Wednesday 18} & {\tiny 4.8} & {\tiny 4.86} & {\tiny 5.005} &
{\tiny 5.1} & {\tiny 4.93} & {\tiny 5.22} & {\tiny 5.335} & {\tiny 5.38} &
{\tiny 5.395} & {\tiny 5.445} & {\tiny 5.64}\\\hline
{\tiny Thursday 19} & {\tiny 4.77} & {\tiny 4.83} & {\tiny 4.97} &
{\tiny 5.06} & {\tiny 4.93} & {\tiny 5.165} & {\tiny 5.27} & {\tiny 5.34} &
{\tiny 5.34} & {\tiny 5.39} & {\tiny 5.585}\\\hline
\end{tabular}
\caption{Quoted rates for the eleven Swedish Government Bonds used in the
calculation of the forward rate curves corresponding to Figs.(\ref{Graph2}%
-\ref{Graph3}). All days are from a two-week period of July 2001. The
description of each of the bonds is given in Table \ref{Table2}. The prices of
the bonds are obtained from this rates and the data of Table \ref{Table2}
using Eq.(\ref{price}). }%
\label{Table1}%
\end{table}

The price of the bond $j$ is calculated using the formula%
\begin{equation}
P_{j}=\sum\limits_{i=1}^{n_{j}}\frac{\left(  c_{j}/100+\delta_{i,n_{j}%
}\right)  N_{j}}{\left(  1+\frac{r_{j}}{100}\right)  ^{\Delta T_{i}}},
\label{price}%
\end{equation}
where $N_{j}$ is the nominal amount (SEK 40 millions for all bonds in Table
\ref{Table1}), $c_{j}$ is the coupon rate (given in Table \ref{Table2}),
$n_{j}$ is the total number of remaining coupons (each paid at time $T_{i}%
^{j}$), $r_{j}$ is the quoted rate given in Table \ref{Table1} and $\Delta
T_{i}$ is a time difference between $T_{i}^{j}$ and the settlement day. This
time difference is calculated according to the ISMA 30E/360 convention defined
as follows: given two dates $\left(  d_{1},m_{1},y_{1}\right)  $ and $\left(
d_{2},m_{2},y_{2}\right)  $, their ISMA 30E/360 time difference $\Delta T$ is
given by%
\begin{align*}
if~~d_{1}  &  =31~~set~~d_{1}~~to~~30,\\
if~~d_{2}  &  =31~~set~~d_{2}~~to~~30,\\
\Delta T  &  =y_{2}-y_{1}+\frac{30\left(  m_{2}-m_{1}\right)  +d_{2}-d_{1}%
}{360},
\end{align*}
\begin{table}[h]
\centering
\begin{tabular}
[c]{|c||c|c|}\hline
{\tiny Bond} & $\overset{\text{{\tiny Maturity}}}{\text{{\tiny (dd/mm/yyyy)}}%
}$ & $\overset{\text{{\tiny annual}}}{\text{{\tiny coupon (\%)}}}%
$\\\hline\hline
{\tiny SO 1033} & {\tiny 05/05/2003} & {\tiny 10.25}\\\hline
{\tiny SO 1042} & {\tiny 15/01/2004} & {\tiny 5}\\\hline
{\tiny SO 1035} & {\tiny 09/02/2005} & {\tiny 6}\\\hline
{\tiny SO 1044} & {\tiny 20/04/2006} & {\tiny 3.5}\\\hline
{\tiny SO 1038} & {\tiny 25/10/2006} & {\tiny 6.5}\\\hline
{\tiny SO 1037} & {\tiny 15/08/2007} & {\tiny 8}\\\hline
{\tiny SO 1040} & {\tiny 05/05/2008} & {\tiny 6.5}\\\hline
{\tiny SO 1043} & {\tiny 28/01/2009} & {\tiny 5}\\\hline
{\tiny SO 1034} & {\tiny 20/04/2009} & {\tiny 9}\\\hline
{\tiny SO 1045} & {\tiny 15/03/2011} & {\tiny 5.25}\\\hline
{\tiny SO 1041} & {\tiny 05/05/2014} & {\tiny 6.75}\\\hline
\end{tabular}
\caption{Maturities and coupon rates for the 11 Swedish bonds quoted in Table
\ref{Table1}. Coupons are paid annually with the last coupon paid at maturity
together with a nominal of SEK 40 millions.}%
\label{Table2}%
\end{table}

\end{document}